\documentclass{article}
\usepackage{amsthm}
\usepackage{amssymb}
\usepackage{amsxtra}

\newtheorem{thm}{Theorem}[section]
\newtheorem{lem}{Lemma}[section]

\begin{document}

\title{On the second moment of $S(T)$ in the theory of the Riemann zeta
function}
\author{Tsz Ho Chan}
\maketitle
\begin{abstract}
We assume the Riemann Hypothesis and an quantitative form of the
Twin Prime Conjecture, and obtain an asymptotic formula for the
second moment of $S(T)$ with better error term.
\end{abstract}

\section{Introduction}
Let $\rho = \beta + i\gamma$ be the zeros of the Riemann zeta
function $\zeta(s)$. For $T \not= \gamma$,
\begin{equation*}
S(T)={1 \over \pi} \mbox{arg} \zeta({1 \over 2} + iT),
\end{equation*}
where the argument is obtained by continuous variation along the
horizontal line $\sigma+iT$ starting with the value zero at
$\infty+ iT$. For $T = \gamma$, we define
\begin{equation*}
S(T) = \lim_{\epsilon \rightarrow 0} {1 \over 2} \{ S(T+\epsilon)
+ S(T-\epsilon) \}.
\end{equation*}
In [\ref{G}], Goldston proved that under the Riemann Hypothesis,
\begin{eqnarray}
\label{1}
\int_{0}^{T} |S(t)|^2 dt &=& {T \over 2\pi^2}
\log{\log{T}} + {T \over 2\pi^2} \Bigl[\int_{1}^{\infty}
{F(\alpha, T) \over \alpha^2} + C_0 \nonumber \\
& &+ \sum_{m=2}^{\infty} \sum_{p} \Bigl({-1 \over m} + {1 \over
m^2}\Bigr) {1 \over p^m} \Bigr] + o(T)\, .
\end{eqnarray}
Here $C_0$ is Euler's constant, and
\begin{equation}
\label{pair} F(\alpha) = F(\alpha, T) = \Bigl({T \over 2\pi}
\log{T}\Bigr)^{-1} \sum_{0 < \gamma, \gamma' \leq T} T^{i \alpha
(\gamma-\gamma')} w(\gamma-\gamma')
\end{equation} with $w(u) = 4/(4+u^2)$ is Montgomery's pair
correlation function. Here and throughout this paper, $p$ will
denote a prime and sums over $p$ are over all primes.

In [\ref{C1}] and [\ref{C2}], the author assumed the Riemann
Hypthesis and the following quantitative form of Twin Prime
Conjecture: For any $\epsilon > 0$,
\begin{equation*}
\sum_{n=1}^{N} \Lambda(n) \Lambda(n+d) = {\mathfrak S}(d) N +
O(N^{1/2+\epsilon})
\end{equation*}
uniformly in $|d| \leq N$. $\Lambda(n)$ is the von Mangoldt lambda
function. ${\mathfrak S}(d) = 2\prod_{p>2} \bigl(1-1/ (p-1)^2
\bigr) \prod_{p|d, p>2} (p-1)/(p-2)$ if $d$ is even, and
${\mathfrak S}(d) = 0$ if $d$ is odd.

\smallskip

He proved that, for any $\epsilon > 0$,
\begin{equation}
\label{2} F(\alpha, T) = \left\{ \begin{array}{ll} \alpha + {1
\over T^{2\alpha}}[\log{T} - 2\log{2\pi} -2] \\
+ O(\alpha T^{\alpha-1}) + O\Bigl({1 \over T^{(1/2-\epsilon)
\alpha} \log{T}} \Bigr), & \raisebox{1ex}{if $0 \leq \alpha \leq 1
-{3\log{\log{T}} \over \log{T}}$,}\\
\raisebox{-1ex}{$\alpha + O({T^{\alpha-1} \over \log{T}})$,} &
\raisebox{-1ex}{if $1-{3\log{\log{T}} \over \log{T}} \leq \alpha
\leq 1$.}
\end{array} \right.
\end{equation}
Using this and more careful calculations, we have
\begin{thm}
\label{main} Assume the Riemann Hypothesis and the above Twin
Prime Conjecture. We have
\begin{eqnarray}
\int_{0}^{T} |S(t)|^2 dt &=& {T \over 2\pi^2} \log{\log{T}} + {T
\over 2\pi^2} \Bigl[\int_{1}^{\infty} {F(\alpha, T) \over
\alpha^2} d\alpha + C_0 \nonumber \\
& & - \sum_{m=2}^{\infty} \sum_{p} \Bigl({1 \over m} - {1 \over
m^2}\Bigr) {1 \over p^m} \Bigr] + O\Bigl(\frac{T}{\log^2{T}}
\Bigr). \nonumber
\end{eqnarray}
\end{thm}
So, there appears to be no $T/\log{T}$ term. However, one has to
be careful about this because we still do not know anything very
precise about
\begin{equation*}
\int_{1}^{\infty} {F(\alpha, T) \over \alpha^2} d\alpha.
\end{equation*}
Goldston [\ref{G}] proved that under the Riemann Hypothesis, for
any $\epsilon > 0$,
\begin{equation*}
\frac{2}{3} - \epsilon < \int_{1}^{\infty} {F(\alpha, T) \over
\alpha^2} d\alpha < 2.
\end{equation*}
Montgomery [\ref{M}] conjectured that
\begin{equation*}
F(\alpha, T) = 1 + o(1) \mbox{ uniformly for } 1 \leq \alpha \leq
M,
\end{equation*}
for any fixed $M$. This implies
\begin{equation*}
\int_{1}^{\infty} {F(\alpha, T) \over \alpha^2} d\alpha = 1 + o(1)
\end{equation*}
which gives an error term $o(T)$ for our Theorem \ref{main}. In
order to get better error terms, one has to understand $F(\alpha,
T)$ better for larger range of $\alpha$, say $1 \leq \alpha \leq
\log{T}$. This would be a great challenge.
\section{Preparations}
Essentially, our proof follows that of Goldston [\ref{G}]. We
shall recall some of his results.
\begin{lem}
\label{lemma1} Assume the Riemann Hypothesis. For $t \geq 1$, $t
\not= \gamma$, $x \geq 4$, we have
\begin{equation}
\label{S(t)}
\begin{split}
S(t) =& -\frac{1}{\pi} \sum_{n \leq x} \frac{\Lambda(n)}{n^{1/2}}
\frac{\sin{(t\log{n})}}{\log{n}} f\Bigl(\frac{\log{n}}{\log{x}}
\Bigr) \\
&+ \frac{1}{\pi} \sum_{\gamma} \sin{((t-\gamma)\log{x})}
\int_{0}^{\infty} \frac{u}{u^2 + ((t-\gamma)\log{x})^2}
\frac{du}{\sinh{u}} \\
&+ O\Bigl(\frac{x^{1/2}}{t^2 \log{x}}\Bigr) + O\Bigl(\frac{1}{t
\log{x}}\Bigr),
\end{split}
\end{equation}
where
\begin{equation}
\label{f(u)} f(u) = \frac{\pi}{2} u \cot{\Bigl(\frac{\pi}{2} u
\Bigr)},
\end{equation}
and $\Lambda(n) = \log{p}$ if $n=p^m$, for $p$ a prime and $m \geq
1$, and $\Lambda(n) = 0$ otherwise.
\end{lem}

Proof: This is Lemma $1$ in [\ref{G}].

\begin{lem}
\label{lemma2} For $x \geq 4$ and $T \geq 2$,
\begin{eqnarray*}
R &=& \int_{1}^{T} \Big| \frac{1}{\pi} \sum_{\gamma}
\sin{((t-\gamma)\log{x})} \int_{0}^{\infty} \frac{u}{u^2 +
((t-\gamma)\log{x})^2} \frac{du}{\sinh{u}} \Big|^2 dt \\
&=& \frac{1}{\pi^2 \log{x}} \sum_{\gamma, \gamma'} \hat{k}((\gamma
- \gamma') \log{x}) + O(\log^3{T}),
\end{eqnarray*}
where
\[ k(u) = \left\{ \begin{array}{ll}
\Bigl(\frac{1}{2u} - \frac{\pi^2}{2} \cot{(\pi^2 u)} \Bigr)^2, &
\mbox{ if } |u| \leq \frac{1}{2\pi} \\
\frac{1}{4u^2}, & \mbox{ if } |u| > \frac{1}{2\pi}.
\end{array} \right. \]
\end{lem}

Proof: This is Lemma $2$ in [\ref{G}].

\begin{lem}
\label{lemma3}
\begin{equation*}
\begin{split}
k'(0) = 0, & \mbox{ } k''(0) = \frac{\pi^8}{18}, \\
k'(\frac{1^+}{2\pi}) = -4 \pi^3, &\mbox{ } k'(\frac{1^-}{2\pi}) =
-4\pi^3 + \pi^5, \\
k''(\frac{1^+}{2\pi}) = 24 \pi^4, & \mbox{ } k''(\frac{1^-}{2\pi})
= \frac{\pi^8}{2} - 4 \pi^6 + 24 \pi^4.
\end{split}
\end{equation*}
\end{lem}

Proof: Straightforward calculations.

\begin{lem}
\label{lemma4}
\begin{equation*}
\hat{k}(y) = -\frac{1}{(2\pi y)^2} \int_{-\infty}^{\infty} k''(u)
e(-uy) du + \frac{\pi^3}{2 y^2} \cos{y},
\end{equation*}
where $\hat{f}$ is the Fourier transform of $f$,
\begin{equation*}
\hat{f}(y) = \int_{-\infty}^{\infty} f(u) e(-uy) du, \mbox{ } e(u)
= e^{2\pi i u}.
\end{equation*}
\end{lem}

Proof: Note that $k(u) \ll \mbox{min}(1, 1/u^2)$, $k'(u) \ll
\mbox{min}(1, 1/u^3)$ and $k''(u) \ll \mbox{min}(1, 1/u^4)$ except
at $u = 1/2\pi$. Integrating by parts twice,
\begin{eqnarray*}
\hat{k}(y) &=& \int_{-\infty}^{\infty} k(u) e(-u y) du \\
&=& \frac{-1}{2 \pi i y} \int_{-\infty}^{\infty} k(u) d e(-u y) =
\frac{1}{2 \pi i y} \int_{-\infty}^{\infty} e(-u y) k'(u) du \\
&=& \frac{-1}{(2\pi i y)^2} \int_{-\infty}^{\infty} k'(u) d e(-u
y) = \frac{1}{(2\pi i y)^2} \int_{-\infty}^{\infty} e(-u y) d
k'(u) \\
&=& \frac{1}{(2\pi i y)^2} \int_{-\infty}^{\infty} k''(u) e(-u y)
du \\
& & + \Bigl(k'(\frac{1^+}{2\pi}) - k'(\frac{1^-}{2\pi})\Bigr)
e(\frac{y}{2\pi}) + \Bigl(k'(-\frac{1^+}{2\pi}) -
k'(-\frac{1^-}{2\pi})\Bigr) e(-\frac{y}{2\pi})\\
&=& \frac{-1}{(2\pi y)^2} \int_{-\infty}^{\infty} k''(u) e(-u y)
du + \frac{\pi^3}{2 y^2} \cos{y}
\end{eqnarray*}
by Lemma \ref{lemma3} and the fact that $k(u)$ is even.

\smallskip

Our key improvement is the following
\begin{lem}
\label{lemma5} Let $x = T^\beta$. For any $\beta > 0$,
\begin{equation*}
\begin{split}
& \sum_{0 < \gamma, \gamma' \leq T} \hat{k}((\gamma-\gamma')
\log{x}) \frac{(\gamma - \gamma')^2}{4 + (\gamma - \gamma')^2} \\
=& \frac{\pi^2 T}{16 \log{T}} \frac{F(\beta)}{\beta^2} -
\frac{T}{64 \pi^4 \log{T} \beta^3} \int_{-\infty}^{\infty}
F(\alpha) k'' \Bigl(\frac{\alpha}{2 \pi \beta}\Bigr) d\alpha,
\end{split}
\end{equation*}
where $F(\alpha)$ is as defined in (\ref{pair}).
\end{lem}

Proof: By Lemma \ref{lemma4}, the above sum
\begin{eqnarray*}
&=& \frac{\pi^3}{2 (\log{x})^2} \sum_{0 < \gamma, \gamma' \leq T}
\cos{((\gamma-\gamma') \log{x})} \frac{1}{4 + (\gamma-\gamma')^2}
\\
& & - \frac{1}{(2\pi \log{x})^2} \sum_{0 < \gamma, \gamma' \leq T}
\int_{-\infty}^{\infty} k''(u) e(-u(\gamma - \gamma') \log{x}) du
\frac{1}{4+(\gamma-\gamma')^2} \\
&=& \frac{\pi^3}{8 (\log{x})^2} \sum_{0 < \gamma, \gamma' \leq T}
x^{i(\gamma-\gamma')} w(\gamma-\gamma') \\
& &- \frac{1}{(4 \pi \log{x})^2} \int_{-\infty}^{\infty} \Bigl(
\sum_{0 < \gamma, \gamma' \leq T} e(-u(\gamma-\gamma') \log{x})
w(\gamma -\gamma')
\Bigr) k''(u) du \\
&=& \frac{\pi^3 \frac{T}{2\pi} \log{T}}{8 (\log{x})^2} F(\beta) -
\frac{\frac{T}{2\pi}\log{T}}{(4 \pi \log{x})^2}
\int_{-\infty}^{\infty} F(2\pi u \beta) k''(u) du
\end{eqnarray*}
which gives the lemma by substituting $\alpha = 2 \pi u \beta$ in
the integral.

Using Lemma \ref{lemma2} and Lemma \ref{lemma5}, we have the
following improvement of Lemma $3$ in [\ref{G}]:
\begin{lem}
\label{lemma6} Let $x = T^\beta$. Then for $\beta > 0$,
\begin{eqnarray*}
R &=& \frac{T}{(2 \pi^2 \beta)^2} \int_{-\infty}^{\infty}
F(\alpha) k\Bigl(\frac{\alpha}{2\pi \beta}\Bigr) d\alpha +
\frac{T}{16 \log^2{T}} \frac{F(\beta)}{\beta^3} \\
& & - \frac{T}{64 \pi^6 \beta^4 \log^2{T}} \int_{-\infty}^{\infty}
F(\alpha) k''\Bigl(\frac{\alpha}{2\pi \beta}\Bigr) d\alpha +
O(\log^3{T}).
\end{eqnarray*}
\end{lem}

Proof: We simply note that
\begin{eqnarray*}
\sum_{0< \gamma, \gamma' \leq T} \hat{k}((\gamma-\gamma') \log{x})
&=& \sum_{0< \gamma, \gamma' \leq T} \hat{k}((\gamma-\gamma')
\log{x}) w(\gamma-\gamma') \\
& & + \sum_{0< \gamma, \gamma' \leq T} \hat{k}((\gamma-\gamma')
\log{x}) \frac{(\gamma-\gamma')^2}{4 + (\gamma-\gamma')^2}
\end{eqnarray*}
Using Lemma \ref{lemma2} and Lemma \ref{lemma5}, we get the second
and third terms. Meanwhile,
\begin{eqnarray*}
& & \sum_{0< \gamma, \gamma' \leq T} \hat{k}((\gamma-\gamma')
\log{x}) w(\gamma-\gamma') \\
&=& \int_{-\infty}^{\infty} k(u) \sum_{0 < \gamma, \gamma' \leq T}
e(-u(\gamma-\gamma') \log{x}) w(\gamma-\gamma') du \\
&=& \frac{T}{2\pi} \log{T} \int_{-\infty}^{\infty} F(2 \pi u
\beta) k(u) du \\
&=& \frac{T \log{T}}{(2\pi)^2 \beta} \int_{-\infty}^{\infty}
F(\alpha) k\Bigl(\frac{\alpha}{2 \pi \beta}\Bigr) d\alpha.
\end{eqnarray*}
which accounts for the first term.

\begin{lem}
\label{lemma7} For any $\beta > 0$,
\begin{equation*}
\int_{0}^{\beta} T^{-2\alpha} k''\Bigl(\frac{\alpha}{2\pi \beta}
\Bigr) d\alpha = 16 \pi^2 \beta^2 (\log{T})^2 \int_{0}^{\beta}
T^{-2\alpha} k\Bigl(\frac{\alpha}{2\pi \beta} \Bigr) d\alpha .
\end{equation*}
\end{lem}

Proof: Integrating by parts twice.

\begin{lem}
\label{lemma8} Assume the Riemann Hypothesis and Twin Prime
Conjecture. For any $\epsilon > 0$ and $0 < \beta < 1$,
\begin{eqnarray*}
\int_{-\infty}^{\infty} F(\alpha) k\Bigl(\frac{\alpha}{2 \pi
\beta} \Bigr) d\alpha &=& 2 \pi^2 \beta^2 \Bigl[1 -
\frac{\pi^2}{8} + \log{\frac{\pi}{2}} + \int_{1}^{\infty}
\frac{F(\alpha)}{\alpha^2} d\alpha - \log{\beta} \Bigr]\\
& & + 2 (\log{T} + C) \int_{0}^{\beta} T^{-2\alpha}
k\Bigl(\frac{\alpha}{2\pi \beta}\Bigr) d\alpha \\
& &+ O\Bigl(\frac{1}{\beta^2 \log^4{T}}\Bigr) +
O\Bigl(\frac{\log{T}}{T^{(1/2-\epsilon)\beta}}\Bigr) +
O\Bigl(\frac{\beta^2}{\log^2{T}}\Bigr)
\end{eqnarray*}
with $C = -2\log{2\pi} - 2$.
\end{lem}

Proof: Let $\epsilon_T = 3\log{\log{T}} / \log{T}$. Since $F$ and
$k$ are even,
\begin{equation*}
\begin{split}
\int_{-\infty}^{\infty} F(\alpha) k\Bigl(\frac{\alpha}{2 \pi
\beta} \Bigr) d\alpha =& 2 \Bigl( \int_{0}^{\beta} +
\int_{\beta}^{1-\epsilon_T} + \int_{1-\epsilon_T}^{1} +
\int_{1}^{\infty} \Bigr) F(\alpha)
k\Bigl(\frac{\alpha}{2 \pi \beta} \Bigr) d\alpha \\
=& 2 (I_1 + I_2 + I_3 + I_4)
\end{split}
\end{equation*}
From (\ref{2}), with $C = -2\log{2\pi} - 2$,
\begin{equation*}
\begin{split}
I_1 =& \int_{0}^{\beta} \Bigl[\alpha + \frac{log{T} + C}{T^{2
\alpha}}\Bigr] \Bigl[\frac{\pi \beta}{\alpha} - \frac{\pi^2}{2}
\cot{(\frac{\pi \alpha}{2 \beta})} \Bigr]^2 d\alpha \\
&+ O\Bigl(\int_{0}^{\beta} \alpha T^{\alpha-1}
\frac{\alpha^2}{\beta^2} d\alpha\Bigr) + O\Bigl(\int_{0}^{\beta}
\frac{T^{-(1/2-\epsilon) \alpha}}{\log{T}}
\frac{\alpha^2}{\beta^2} d\alpha\Bigr) \\
=& \int_{0}^{\beta} \alpha \Bigl[\frac{\pi \beta}{\alpha} -
\frac{\pi^2}{2} \cot{(\frac{\pi \alpha}{2 \beta})} \Bigr]^2
d\alpha + (\log{T} + C) \int_{0}^{\beta} T^{-2\alpha}
k\Bigl(\frac{\alpha}{2\pi \beta}\Bigr) d\alpha\\
&+ O\Bigl(\frac{T^{\beta-1}}{\beta^2 \log^4{T}}\Bigr) +
O\Bigl(\frac{1}{\beta^2 \log^4{T}}\Bigr)
\end{split}
\end{equation*}
because $\cot(x) = 1/x + O(x)$ when $0 \leq x \leq \pi/2$. The
first integral is elementary to evaluate.
\begin{equation*}
I_1 = \pi^2 \beta^2 \Bigl[1 - \frac{\pi^2}{8} +
\log{\frac{\pi}{2}}\Bigr] + (\log{T} + C) \int_{0}^{\beta}
T^{-2\alpha} k\Bigl(\frac{\alpha}{2\pi \beta}\Bigr) d\alpha +
O\Bigl(\frac{1}{\beta^2 \log^4{T}}\Bigr).
\end{equation*}
By (\ref{2}) again, we have,
\begin{equation*}
\begin{split}
I_2 =& \int_{\beta}^{1-\epsilon_T} \Bigl[\alpha + O(\alpha
T^{\alpha-1}) + O\Bigl(\frac{\log{T}}{T^{(1/2-\epsilon)\alpha}}
\Bigr)\Bigr] \Bigl(\frac{\pi \beta}{\alpha}\Bigr)^2 d\alpha \\
=& \pi^2 \beta^2 [\log{(1-\epsilon_T)} - \log{\beta}] +
O\Bigl(\frac{1}{\log^4{T}}\Bigr) +
O\Bigl(\frac{\log{T}}{T^{(1/2-\epsilon)\beta}}
\Bigr) \\
I_3 =& \int_{1-\epsilon_T}^{1} \Bigl[\alpha +
O\Bigl(\frac{T^{\alpha-1}}{\log{T}}\Bigr) \Bigr] \Bigl(\frac{\pi
\beta}{\alpha}\Bigr)^2 d\alpha \\
=& - \pi^2 \beta^2 \log{(1-\epsilon_T)} +
O\Bigl(\frac{\beta^2}{\log{T}} \int_{1-\epsilon_T}^{1}
T^{\alpha-1} d\alpha \Bigr) \\
=& - \pi^2 \beta^2 \log{(1-\epsilon_T)} +
O\Bigl(\frac{\beta^2}{\log^2{T}}\Bigr).
\end{split}
\end{equation*}
Finally,
\begin{equation*}
I_4 = \pi^2 \beta^2 \int_{1}^{\infty} \frac{F(\alpha)}{\alpha^2}
d\alpha.
\end{equation*}
Combining the results for $I_1$, $I_2$, $I_3$ and $I_4$, we have
the lemma.
\begin{lem}
\label{lemma9} Assume the Riemann Hypothesis and Twin Prime
Conjecture. For any $\epsilon > 0$ and $0 < \beta < 1$, where
$\beta = \log{x}/\log{T}$,
\begin{equation*}
\begin{split}
\int_{-\infty}^{\infty} F(\alpha) k''\Bigl(\frac{\alpha}{2 \pi
\beta}\Bigr) d\alpha =& 4 \pi^6 \beta^2 - 24 \pi^4 \beta^4 + 48
\pi^4 \beta^4 \int_{1}^{\infty} \frac{F(\alpha)}{\alpha^4} d\alpha
\\
&+ 32 \pi^2 \beta^2 (\log{T})^2 (\log{T} + C) \int_{0}^{\beta}
T^{-2\alpha} k\Bigl(\frac{\alpha}{2 \pi \beta} \Bigr) d\alpha \\
&+ O\Bigl(\frac{1}{\log^2{T}}\Bigr) +
O\Bigl(\frac{\log{T}}{T^{(1/2-\epsilon)\beta}}\Bigr)
\end{split}
\end{equation*}
\end{lem}

Proof: Let $\epsilon_T = 3\log{\log{T}} / \log{T}$. Again, since
$F$ and $k$ are even,
\begin{equation*}
\begin{split}
\int_{-\infty}^{\infty} F(\alpha) k''\Bigl(\frac{\alpha}{2 \pi
\beta}\Bigr) d\alpha =& 2\Bigl(\int_{0}^{\beta} +
\int_{\beta}^{1-\epsilon_T} + \int_{1-\epsilon_T}^{1} +
\int_{1}^{\infty} \Bigr) F(\alpha)
k''\Bigl(\frac{\alpha}{2 \pi \beta}\Bigr) d\alpha \\
=& 2(J_1 + J_2 + J_3 + J_4)
\end{split}
\end{equation*}
By (\ref{2}),
\begin{equation*}
\begin{split}
J_1 =& \int_{0}^{\beta} \alpha k''\Bigl(\frac{\alpha}{2 \pi
\beta}\Bigr) d\alpha + 16 \pi^2 \beta^2 (\log{T})^2 (\log{T} + C)
\int_{0}^{\beta} T^{-2\alpha} k\Bigl(\frac{\alpha}{2 \pi \beta}
\Bigr) d\alpha \\
&+ O\Bigl(\int_{0}^{\beta} \alpha T^{\alpha-1} d\alpha\Bigr) +
O\Bigl(\int_{0}^{\beta} \frac{T^{-(1/2-\epsilon)\alpha}}{\log{T}}
d\alpha\Bigr) \\
=& 4\pi^2 \beta^2 \int_{0}^{1/2\pi} u k''(u) du + 16 \pi^2 \beta^2
(\log{T})^2 (\log{T} + C) \int_{0}^{\beta} T^{-2\alpha}
k\Bigl(\frac{\alpha}{2 \pi \beta} \Bigr) d\alpha \\
&+ O\Bigl(\frac{T^{\beta-1}}{\log^2{T}}\Bigr) +
O\Bigl(\frac{1}{\log^2{T}}\Bigr)
\end{split}
\end{equation*}
since $k''(x) \ll 1$ when $0 \leq x \leq 1/2\pi$. Using
integration by parts twice and Lemma \ref{lemma3} to compute the
first integral, we have
\begin{equation*}
\begin{split}
J_1 =& 2\pi^4 (\pi^2 - 6) \beta^2 + 16 \pi^2 \beta^2 (\log{T})^2
(\log{T} + C) \int_{0}^{\beta} T^{-2\alpha} k\Bigl(\frac{\alpha}{2
\pi \beta} \Bigr) d\alpha \\
&+ O\Bigl(\frac{1}{\log^2{T}}\Bigr).
\end{split}
\end{equation*}
By (\ref{2}) again, we have
\begin{equation*}
\begin{split}
J_2 =& \int_{\beta}^{1-\epsilon_T} \Bigl[\alpha + O(\alpha
T^{\alpha-1}) + O\Bigl(\frac{\log{T}}{T^{(1/2-\epsilon)\alpha}}
\Bigr)\Bigr] \Bigl(\frac{24 \pi^4 \beta^4}{\alpha^4}\Bigr) d\alpha \\
=& 12 \pi^4 \beta^4 \Bigl[\frac{1}{\beta^2} -
\frac{1}{(1-\epsilon_T)^2}\Bigr] +
O\Bigl(\frac{1}{\log^4{T}}\Bigr) +
O\Bigl(\frac{\log{T}}{T^{(1/2-\epsilon)\beta}} \Bigr) \\
J_3 =& \int_{1-\epsilon_T}^{1} \Bigl[\alpha +
O\Bigl(\frac{T^{\alpha-1}}{\log{T}}\Bigr) \Bigr] \Bigl(\frac{24
\pi^4 \beta^4}{\alpha^4}\Bigr) d\alpha \\
=& 12 \pi^4 \beta^4 \Bigl[\frac{1}{(1-\epsilon_T)^2} - 1 \Bigr] +
O\Bigl(\frac{\beta^4}{\log{T}} \int_{1-\epsilon_T}^{1}
T^{\alpha-1} d\alpha \Bigr) \\
=& 12 \pi^4 \beta^4 \Bigl[\frac{1}{(1-\epsilon_T)^2} - 1 \Bigr] +
O\Bigl(\frac{\beta^4}{\log^2{T}}\Bigr).
\end{split}
\end{equation*}
Finally,
\begin{equation*}
J_4 = 24 \pi^4 \beta^4 \int_{1}^{\infty}
\frac{F(\alpha)}{\alpha^4} d\alpha.
\end{equation*}
Combining the results for $J_1$, $J_2$, $J_3$ and $J_4$, we have
the lemma.

\smallskip

Combining Lemma \ref{lemma6}, Lemma \ref{lemma8} and Lemma
\ref{lemma9}, we have
\begin{lem}
\label{lemma10} Assume the Riemann Hypothesis and Twin Prime
Conjecture, for fixed $0 < \beta < 1$, where $\beta =
\log{x}/\log{T}$,
\begin{equation*}
\begin{split}
R =& \frac{T}{2 \pi^2} \Bigl[1 - \frac{\pi^2}{8} +
\log{\frac{\pi}{2}} + \int_{0}^{\infty} \frac{F(\alpha)}{\alpha^2}
d\alpha - \log{\beta}\Bigr] + \frac{3 T}{8\pi^2 \log^2{T}} \\
&- \frac{3 T}{4\pi^2 \log^2{T}} \int_{1}^{\infty}
\frac{F(\alpha)}{\alpha^4} d\alpha + O\Bigl(\frac{T}{\log^2{T}}
\Bigr) + O\Bigl(\frac{T}{\beta^4 \log^4{T}}\Bigr)
\end{split}
\end{equation*}
\end{lem}

Note: This is a more precise version of Lemma $4$ of [\ref{G}].
Also, we keep some of the $T / \log^2{T}$ terms explicit because
one can actually make the $O(T/\log^2{T})$ error term $=C_1
T/\log^2{T} + O(T \log{\log{T}}/ \log^3{T})$ for some constant
$C_1$ by using Theorem 1.1 in [\ref{C2}].

Following [\ref{G}], we need to compute the mean value of the
Dirichlet series in Lemma \ref{lemma1}, and the cross term
obtained from multiplying $S(t)$ with this series. Let
\begin{equation*}
G(T) = \int_{1}^{T} \Big| \frac{1}{\pi} \sum_{n \leq x}
\frac{\Lambda(n)}{n^{1/2}} \frac{\sin{(t\log{n})}}{\log{n}}
f\Bigl(\frac{\log{n}}{\log{x}}\Bigr) \Big|^2 dt
\end{equation*}
and
\begin{equation*}
H(T) = \frac{2}{\pi} \int_{1}^{T} S(t) \sum_{n \leq x}
\frac{\Lambda(n)}{n^{1/2}} \frac{\sin{(t\log{n})}}{\log{n}}
f\Bigl(\frac{\log{n}}{\log{x}}\Bigr) dt,
\end{equation*}
where $f$ is defined as in (\ref{f(u)}). We need a lemma.
\begin{lem}
\label{lemma11} For $C \geq 2$ and $k \geq 1$,
\begin{equation*}
\sum_{n=1}^{\infty} \frac{n^k}{C^n} \ll_k \frac{1}{C}.
\end{equation*}
\end{lem}

Proof: First, we note that $u^k C^{-u}$ is decreasing when $u >
\frac{k}{\log{C}}$. So,
\begin{eqnarray*}
\sum_{n=1}^{\infty} \frac{n^k}{C^n} &=& \sum_{n=1}^{k/\log{C}}
\frac{n^k}{C^n} + \sum_{n > k/\log{C}} \frac{n^k}{C^n} \\
&\leq& \Bigl(\frac{k}{\log{2}}\Bigr)^k \frac{1}{C-1} +
\int_{1}^{\infty} u^k C^{-u} du \\
&\ll_k& \frac{1}{C} + \frac{1}{\log{C}} C^{-1}
\end{eqnarray*}
by integration by parts. This gives the lemma.

From p.165-166 of [\ref{G}], we have, assuming the Riemann
Hypothesis,
\begin{eqnarray}
\label{G(T)} G(T) &=& \frac{T}{2\pi^2} \sum_{n \leq
x}\frac{\Lambda^2(n)}{n
\log^2{n}} f^2\Bigl(\frac{\log{n}}{\log{x}}\Bigr) + O(x^2), \\
\label{H(T)} H(T) &=& -\frac{T}{\pi^2} \sum_{n \leq x}
\frac{\Lambda^2(n)}{n \log^2{n}}
f\Bigl(\frac{\log{n}}{\log{x}}\Bigr) + O(x^{2+\epsilon})
\end{eqnarray}
for any $\epsilon > 0$. Adding (\ref{G(T)}) and (\ref{H(T)}),  we
have
\begin{equation*}
\begin{split}
G(T) + H(T) =& \frac{T}{2\pi^2} \Bigl[ \sum_{p \leq x} \frac{1}{p}
f^2\Bigl(\frac{\log{p}}{\log{x}}\Bigr) - 2 \sum_{p \leq x}
\frac{1}{p} f\Bigl(\frac{\log{p}}{\log{x}}\Bigr) -
\sum_{m=2}^{\infty} \sum_{p^m \leq x} \frac{1}{m^2 p^m}\\
&+ \sum_{m=2}^{\infty} \sum_{p^m \leq x} \frac{1}{m^2 p^m}
\Bigl(f\Bigl(\frac{m\log{p}}{\log{x}}\Bigr) -1 \Bigr)^2 \Bigr] +
O(x^{2+\epsilon}) \\
=& \frac{T}{2\pi^2} [S_1 - 2S_2 - S_3 + S_4] + O(x^{2+\epsilon}).
\end{split}
\end{equation*}
\begin{eqnarray*}
S_3 &=& \sum_{m=2}^{\infty} \sum_{p} \frac{1}{m^2 p^m} +
O\Bigl(\sum_{m=2}^{\infty} \frac{1}{m^2} \sum_{n \geq x^{1/m}}
\frac{1}{n^m}\Bigr) \\
&=& \sum_{m=2}^{\infty} \sum_{p} \frac{1}{m^2 p^m} +
O\Bigl(\sum_{m=2}^{\infty} \frac{1}{m^2 (m-1) x^{1-1/m}}\Bigr) \\
&=& \sum_{m=2}^{\infty} \sum_{p} \frac{1}{m^2 p^m} +
O\Bigl(\frac{1}{x^{1/2}}\Bigr).
\end{eqnarray*}
By Taylor's expansion of $\tan{x}$, we have $f(u) = 1 + O(u^2)$
when $0 \leq u \leq 1$. Thus,
\begin{eqnarray*}
S_4 &\ll& \frac{1}{\log^4{x}} \sum_{m=2}^{\infty} m^2 \sum_{p^m
\leq x} \frac{\log^4{p}}{p^m} \\
&\ll& \frac{1}{\log^4{x}} \sum_{m=2}^{\infty} m^2
\sum_{i=1}^{\infty} \sum_{2^i \leq p \leq 2^{i+1}} \frac{i^4}{2^{m
i}} \\
&\ll& \frac{1}{\log^4{x}} \sum_{m=2}^{\infty} m^2
\sum_{i=1}^{\infty} \frac{i^4}{(2^{m-1})^i} \\
&\ll& \frac{1}{\log^4{x}} \sum_{m=2}^{\infty} \frac{m^2}{2^{m-1}}
\ll \frac{1}{\log^4{x}}
\end{eqnarray*}
by using Lemma \ref{lemma11} twice.

We now define
\begin{equation*}
T(u) = \sum_{2 \leq p \leq u} \frac{1}{p},
\end{equation*}
and have
\begin{equation*}
T(u) = \log{\log{u}} + C_0 + \sum_{p} \Bigl(\log{\Bigl(1 -
\frac{1}{p}\Bigr)} + \frac{1}{p}\Bigr) + r(u),
\end{equation*}
where $r(u) \ll \log{u}/\sqrt{u}$ under the Riemann Hypothesis,
and the sum over primes on the right is equal to
\begin{equation*}
-\sum_{m=2}^{\infty} \sum_{p} \frac{1}{m p^m}.
\end{equation*}
Then,
\begin{eqnarray*}
S_1 &=& \int_{2}^{x} f^2\Bigl(\frac{\log{u}}{\log{x}} \Bigr) d
T(u) \\
&=& \int_{2}^{x} f^2\Bigl(\frac{\log{u}}{\log{x}} \Bigr)
\frac{du}{u\log{u}} + \int_{2}^{x}
f^2\Bigl(\frac{\log{u}}{\log{x}}\Bigr) d r(u) = I_1 + I_2.
\end{eqnarray*}
Similarly,
\begin{equation*}
S_2 = \int_{2}^{x} f\Bigl(\frac{\log{u}}{\log{x}} \Bigr)
\frac{du}{u\log{u}} + \int_{2}^{x}
f\Bigl(\frac{\log{u}}{\log{x}}\Bigr) d r(u) = J_1 + J_2.
\end{equation*}
Thus,
\begin{equation*}
S_1 - 2S_2 = I_1 - 2J_1 + (I_2 - 2J_2).
\end{equation*}
By integration by parts,
\begin{eqnarray*}
I_2-2J_2 &=& -r(2^-) \Bigl[f^2\Bigl(\frac{\log{2}}{\log{x}}\Bigr)
- 2f\Bigl(\frac{\log{2}}{\log{x}}\Bigr)\Bigr] - \int_{2}^{x} r(u)
\\
& &\Bigl[2 f\Bigl(\frac{\log{u}}{\log{x}}\Bigr) f'\Bigl(
\frac{\log{u}}{\log{x}}\Bigr) - f'\Bigl(
\frac{\log{u}}{\log{x}}\Bigr) \Bigr] \frac{du}{u\log{x}} \\
&=& r(2^-) - r(2^-) \Bigl[f\Bigl(\frac{\log{2}}{\log{x}}\Bigr) -
1\Bigr]^2 \\
& & -\frac{2}{\log{x}} \int_{2}^{x} r(u) f'\Bigl(
\frac{\log{u}}{\log{x}}\Bigr) \Bigl[f\Bigl(
\frac{\log{u}}{\log{x}}\Bigr) - 1\Bigr] \frac{du}{u} \\
&=& -\log{\log{2}} - C_0 + \sum_{m=2}^{\infty} \sum_{p} \frac{1}{m
p^m} + O\Bigl(\frac{1}{\log^4{x}}\Bigr) \\
& & + O\Bigl(\frac{1}{\log{x}} \int_{2}^{x}
\frac{\log{u}}{\sqrt{u}} \Bigl(\frac{\log{u}}{\log{x}}\Bigr)
\Bigl(\frac{\log{u}}{\log{x}}\Bigr)^2 \frac{du}{u} \Bigr)\\
&=& -\log{\log{2}} - C_0 + \sum_{m=2}^{\infty} \sum_{p} \frac{1}{m
p^m} + O\Bigl(\frac{1}{\log^4{x}}\Bigr)
\end{eqnarray*}
because $f(u) = 1 + O(u^2)$ and $f'(u) \ll u$ when $0 \leq u \leq
1$. The integrals in $I_1$ and $J_1$ are elementary to evaluate.
Using $u\cot{u} = 1 - u^2/3 + O(u^4)$ and $\sin{u} = u - u^3/6 +
O(u^5)$, one has
\begin{eqnarray*}
I_1 &=& \log{\log{x}} - \log{\log{2}} - \frac{\pi^2}{8} + 1 -
\log{\frac{\pi}{2}} + \frac{\pi^2 (\log{2})^2}{12 \log^2{x}}
+ O\Bigl(\frac{1}{\log^4{x}}\Bigr), \\
J_1 &=& \log{\log{x}} - \log{\log{2}} - \log{\frac{\pi}{2}} +
\frac{\pi^2 (\log{2})^2}{24 \log^2{x}} +
O\Bigl(\frac{1}{\log^4{x}}\Bigr).
\end{eqnarray*}
Hence,
\begin{equation*}
S_1 - 2S_2 = -\log{\log{x}} + \log{\frac{\pi}{2}} -
\frac{\pi^2}{8} +1 -C_0 + \sum_{m=2}^{\infty} \sum_{p} \frac{1}{m
p^m} + O\Bigl(\frac{1}{\log^4{x}}\Bigr).
\end{equation*}
Therefore,
\begin{equation}
\label{G+H}
\begin{split}
G(T) + H(T) =& \frac{T}{2\pi^2} \Bigl[ -\log{\log{x}} +
\log{\frac{\pi}{2}} - \frac{\pi^2}{8} +1 - C_0 \\
&+ \sum_{m=2}^{\infty} \sum_{p}
\Bigl(\frac{1}{m}-\frac{1}{m^2}\Bigr) \frac{1}{p^m} \Bigr] +
O\Bigl(\frac{T}{\log^4{x}}\Bigr).
\end{split}
\end{equation}
\section{Proof of Theorem \ref{main}}
Suppose $x = T^\beta$ and $\beta$ is a fixed positive number less
than $1/2$. We have, by Lemma \ref{lemma1}, that (\ref{S(t)})
holds except on a countable set of points. Hence, on squaring both
sides of (\ref{S(t)}) and integrating from $1$ to $T$,
\begin{equation*}
\int_{1}^{T} (S(t))^2 dt + H(T) + G(T) = R + O(T^{1/2} x^{1/2}),
\end{equation*}
where the error term is obtained by Cauchy-Schwarz inequality
since $R \ll T$. The lower limit of integration may be replaced by
zero since $\int_{0}^{1} (S(t))^2 dt \ll 1$. Then, Lemma
\ref{lemma10} and (\ref{G+H}) give the theorem. The author would
like to thank Professor Daniel Goldston for suggestion and
discussion on this problem.


Tsz Ho Chan\\
Case Western Reserve University\\
Mathematics Department, Yost Hall 220\\
10900 Euclid Avenue\\
Cleveland, OH 44106-7058\\
USA\\
txc50@po.cwru.edu

\end{document}